\documentclass[11pt,leqeq]{article}
\usepackage{amssymb,amsthm}
\usepackage{amsmath}
\usepackage{amscd}
\usepackage{xy}
\xyoption{all}
\usepackage[dvips]{graphicx}

\newtheorem{thm}{Theorem}

\newtheorem{cor}[thm]{Corollary}

\newtheorem{prop}[thm]{Proposition} 

\newtheorem{defn}[thm]{Definition}

\newtheorem{exmp}[thm]{Example}

\begin{document}

\setlength{\baselineskip}{14pt}

\title{Symmetric Boolean Algebras}

\author{Rafael D\'\i az and Mariolys Rivas }

\maketitle


\begin{abstract}
We define Boolean algebras in the linear context and study its
symmetric powers.  We give explicit formulae for products in
symmetric Boolean algebras of various dimensions. We formulate
symmetric forms of the inclusion-exclusion principle.

\end{abstract}


\section*{Introduction}

Fix $k$ a field of characteristic zero. A fundamental fact in
mathematics is the existence of a functor $<\  \ >: Set
\longrightarrow vect$ from the category of
sets into the category $k$-vector spaces.  The functor $<\
\ >:  Set \to vect$ sends a set $x$ into $<x>$ the free $k$-vector
space generated by $x$ and sends a map $f: x\to y$ into the linear
transformation $\widehat{f}: <x> \longrightarrow <y>$ whose value
on $i\in x$
is $f(i)$.\\

It is important to notice that both $Set$ and $vect$ are symmetric
monoidal categories with coproduct and that $<\ \ >$ is a monoidal
functor that respects coproducts.  The monoidal structure on $Set$
is Cartesian product $\times$ and the coproduct is disjoint union
$\sqcup$.  The  monoidal structure on $vect$ is tensor product
$\otimes$ and the coproduct is direct sum $\oplus$. The restricted
functor $<\ \ > : FinSet
\longrightarrow vect$ is such that the dimension $dim(< x>)$ of $< x >$
is the cardinality $|x|$ of $x$ for each $x$ in $FinSet$, the
category of
finite sets. \\

Using $<\ \ >$ we can transform (combinatorial) set theoretical
notions into (finite dimensional) linear algebra notions. For
example the linear analogue of a monoid is an associative algebra
since for any monoid $x$ the vector space $<x>$ carries the
structure of an associative algebra. Similarly, the linear
analogue of a group is a Hopf algebra since $<x>$ carries a
structure of a Hopf algebra for any group $x$.\\

Boolean algebras has been known since 1854 and are a cornerstone
of modern mathematics. Despite its widespread range of
applications we believe the right name for Boolean algebras should
have been Boolean monoids.  For most mathematicians the word
algebra implies a linear structure which is certainly not included
in the traditional definition of Boolean algebras. For the
purposes of this paper we find useful to make the distinction
between Boolean monoids and Boolean algebras. The first goal of
this paper is to uncover the linear analogue of Boolean monoids,
i.e., we answer the question: what algebraic structure does $<x>$
carry for any Boolean monoid $x$? We will see that there
are an infinite number of non-isomorphic Boolean algebras. \\

The second goal of this paper is to study the symmetric powers of
Boolean algebras. We compute the structural constants of such
algebras in various dimension, and show that each symmetric
function can be use to formulate a generalization of the
inclusion-exclusion principle for the symmetric powers of Boolean
algebra. Our final goal is to  propose a logical interpretation
for Boolean algebras and pose some open
problems.\\

\section{Boolean monoids vs Boolean algebras}

We recall the definition of Boolean monoids for definiteness and
reader convenience, so that he or she may contrast it with the
definition of Boolean algebras given below.

\smallskip

\begin{defn}
 A Boolean monoid is a set $B$ together with the data
 \begin{enumerate}
 \item Maps
    $\cup:B\times B\to B$,
    $\cap:B\times B\to B$,
    $c:B\to B$  called union, intersection and complement, respectively.
\item Distinguished elements
$e,t\in B$  called the empty and total element,  respectively.
 \end{enumerate}
This data should satisfy the following identities for  $a,b,c\in
B$
\begin{enumerate}
\item $a\cup b=b\cup a$, $a\cap b=b\cap a$.
\item $a\cup(b\cup c)=(a\cup b)\cup c$, $a\cap(b\cap c)=(a\cap b)\cap c$.
\item $a\cap(b\cup c)=(a\cap b)\cup(a\cap c)$, $a\cup(b\cap c)=(a\cup b)\cap(a\cup c)$.
\item $a\cup e=a$, $a\cap t=a$.
\item $a\cup a^{c}=t$, $a\cap a^{c}=e$.
\item $a \cup (a \cap b) = a,$ $a \cap (a \cup b) = a.$
\end{enumerate}
\end{defn}

To any set $x$ we associate the Boolean monoid  $P(x)=\{a\mid
a\subseteq x\}$ where
\begin{enumerate}
\item $a\cup b=\{i\in x\mid i\in a$ or $i\in b\}$.
\item $a\cap b=\{i\in x\mid i\in a$ and $i\in b\}$.
\item $a^{c}=\{i\in x\mid i\notin a\}$.
\item $t$ is $x$ and $e$ is the empty set $\emptyset$.
\end{enumerate}

Let $[n] = \{ 1,...,n \}$ and $S_n$ be the group of permutations
on $n$ letters. we write $P[n]$ instead of $P([n])$ if no
confusion arises. Examples of the form $P(x)$ are essentially the
unique models of
finite Boolean monoids.\\

\begin{thm}
Every finite Boolean monoid is isomorphic to $P(x)$ for a finite
set $x$.
\end{thm}

\begin{proof}
Let $B$ be a Boolean monoid. Define a partial order $\leq$ on $B$
by letting $a \leq b$ if and only if $a\cap b =a$. Let $x$ be the
set of primitive elements or atoms in $B$, i.e.,
$$x = \{a\in A \ \ | \ \  a \neq e \mbox{ and if } b\leq a \mbox{ then }
b=e\}.$$ The map $f:B\longrightarrow P(x)$ given by $f(b)=\{a \in
X
\mid a\leq b\}$ is an isomorphism between $B$ and $P(x)$.
\end{proof}
The Boolean monoids $P(x)$ are described as follows

\begin{thm}
\begin{itemize}
\item{If $B$ and $C$ are Boolean monoids then $B \times C$ is a Boolean monoid.}
\item{$P(x)$ is isomorphic to ${P[1]}^{|x|}$.}
\end{itemize}
\end{thm}

For a $k$-vector space $V$ we use the symmetry map $S:V\otimes
V\to V\otimes V$ given by $S(x\otimes y)=y\otimes x$ for  $x,y\in
V$. The identity map $I:V\to V$ is given by $I(x)=x$ for $x\in V$.
We are ready to define the linear analogue of the notion of
Boolean monoids.

\begin{defn}
 A Boolean algebra is a $k$-vector space $V$
 together with the  data
\begin{enumerate}
\item Linear maps $\cup:V\otimes V\to V$, $\cap:V\otimes V\to V$,
$c:V\to V$ called union, intersection and complement,
respectively.
\item Linear maps $T:k\to V$,
$E:k\to V$  called the empty map and total map,  respectively.
\item Linear map $\bigtriangleup:V\to V\otimes V$  called coproduct.
\item Linear map $ev:V\to k$  called evaluation.
\end{enumerate}
The axioms below hold

\begin{enumerate}
\item $\cup=\cup\circ S$, $\cap=\cap\circ S$.
\item $\cup\circ(\cup\otimes I)=\cup\circ(I\otimes\cup)$,
$\cap\circ(\cap\otimes I)=\cap\circ(I\otimes\cap)$.
\item$\cap\circ(I\otimes\cup)=\cup\circ(\cap\otimes\cap)\circ(I\otimes
S\otimes I)\circ(\bigtriangleup\otimes I\otimes I)$,\\
$\cup\circ(I\otimes\cap)=\cap\circ(\cup\otimes\cup)\circ(I\otimes
S\otimes I)\circ(\bigtriangleup\otimes I\otimes I)$.
\item $\cup\circ(I\otimes E)=I$, $\cap\circ(I\otimes T)=I$.
\item $\cap\circ(I\otimes c)\circ\bigtriangleup=E\circ ev$, $\cup\circ(I\otimes c)\circ\bigtriangleup=T\circ ev$.
\item $\cap \circ (I \otimes \cup) \circ (\bigtriangleup \otimes I)
= I \otimes ev$,  $\cup \circ (I \otimes \cap) \circ
(\bigtriangleup \otimes I) = I \otimes ev.$
\item $(\bigtriangleup\otimes I)\circ\bigtriangleup=(I\otimes\bigtriangleup)\circ\bigtriangleup$.
\item $S\circ\bigtriangleup=\bigtriangleup$.
\end{enumerate}
\end{defn}

Our next result guarantees the existence of infinitely many models
of Boolean algebras, namely those naturally associated with
Boolean monoids.

\begin{thm}
$<P(x)>$ is a Boolean algebra for any set $x$.
\end{thm}
\begin{proof}
The structural maps are given by
\begin{enumerate}
\item  $(\Sigma_{a \subseteq x} v_{a}a)\cup (\Sigma_{b \subseteq x} w_{b}b)=\Sigma_{a \subseteq x, b \subseteq x} v_{a}w_{b}(a\cup b).$
\item $(\Sigma_{a \subseteq x} v_{a}a)\cap (\Sigma_{b \subseteq x} w_{b}b)=\Sigma_{a \subseteq x, b \subseteq x} v_{a}w_{b}(a\cap b).$
\item $c(\Sigma_{a \subseteq x} v_{a}a)=\Sigma_{a \subseteq x} v_a c(a).$
\item $\bigtriangleup(\Sigma_{a \subseteq x} v_{a}a)=\Sigma_{a \subseteq x} v_{a}(a \otimes a).$
\item $e(s)= s\emptyset,$ for $s \in k$.
\item $t(s) = s x$, for $s \in k$.
\item $ev(\Sigma_{a \subseteq x} v_{a}a)=\Sigma_{a \subseteq x} v_a.$
\end{enumerate}
\end{proof}

Next result characterizes finite dimensional Boolean algebras of
the form $<P(x)>$.

\begin{thm}\label{ba}
\begin{itemize}
\item{If $V$ and $W$ are Boolean algebras then $V \otimes W$ is a Boolean algebra.}
\item{$<P(x)>=<P[1]>^{\otimes |x|}$.}
\end{itemize}
\end{thm}

\begin{proof} The Boolean operations on $V \otimes W$ are
define component wise.
\end{proof}

Theorem \ref{ba} suggests the following open problems.\\

\noindent {\bf Problem 1.}  Is any Boolean algebra isomorphic
to $<B>$ for some Boolean monoid $B$? \\

\noindent {\bf Problem 2.}  Classify all finite dimensional
Boolean algebras. \\

\noindent {\bf Problem 3.} Is any finite dimensional Boolean algebra isomorphic
to $<P[n]>$ for some $n \in \mathbb{N}$? \\

\section{Boolean prop}

In this section we give a scientific explanation for our choice of
axioms for Boolean algebras. We do so by defining the prop in
$vect$ whose algebras are Boolean algebras and showing that this
prop actually comes from a prop in $Set$ whose algebras are
Boolean monoids. Discovering the prop that controls a given family
of algebras is like unveiling its genetic code
\cite{a}, \cite{lane},
\cite{Merk3}, \cite{Merk1},
\cite{Merk2}. Despite the fact that Boolean algebras have been
extensively studied from a myriad of view points its genetic code
has not been study so far. Since the theory of props is not widely
known we provide an overview using a convenient notation for our
purposes. We define props over a symmetric monoidal category
$\mathcal{C}$
\footnote{For technical reasons we assume that objects of $\mathcal{C}$ are
sets and that $\mathcal{C}(1_\mathcal{C},x)=x$ for  $x$ an object
of $\mathcal{C}.$ We also assume that $\mathcal{C}$ admits finite
colimits. }, but the reader should bear in mind that in this work
$\mathcal{C}$ is either $Set$ or $vect$.

\begin{defn}
\begin{itemize}
\item A prop over $\mathcal{C}$ is a symmetric monoidal category $P$
enriched over $\mathcal{C}$ such that 1) $Ob(P)=\mathbb{N}$. 2)
The monoidal structure is addition.

\item Let Prop$_{\mathcal{C}}$ be the category whose objects are
props over $\mathcal{C}.$ Morphisms in Prop$_{\mathcal{C}}$ are
monoidal functors.
\end{itemize}
\end{defn}

By definition each prop $P$ is provided with the following
data

\begin{itemize}
\item For each $n \in \mathbb{N}$, a group morphisms $S_{n}\longrightarrow P(n,n)$
such that the diagram
\[\xymatrix @R=.4in  { S_{n}\times S_{m}
\ar[d]_{}
\ar[r]^{} & S_{n + m} \ar[d]^-{}\\
P(n,n )\otimes_{\mathcal{C}} P(m,m)  \ar[r]_{} & P(n + m, n +
m)}\] commutes. The maps $S_{n}\longrightarrow P(n,n)$ induce a
right  action of $S_n$ on $P(n,m)$ and a left of $S_m$. on
$P(n,m).$

\item Let $\mathbb{B}$ be the category whose objects are finite sets and whose morphisms
are bijections. The actions constructed above are used to define a
functor $P:
\mathbb{B}^{op}
\times \mathbb{B} \longrightarrow \mathcal{C}$ given by
$$P(a,b)=  \mathbb{B}(a,[|a|]) \times_{S_{|b|}} P(|a|,|b|)\times_{S_{|b|}}\mathbb{B}([|b|],b).$$

\item Morphisms $P(n,m)\otimes_{\mathcal{C}} P(m,k)\longrightarrow P(n,k)$ for $n,m,k
\in \mathbb{N}.$

\item  Morphisms $P(n,m)\otimes_{\mathcal{C}} P(k,l)\longrightarrow P(n + k,m +
l)$ for $n,m,k,l \in \mathbb{N}$.

\end{itemize}

In order to define the free prop generated by a functor
$G:\mathbb{B}^{op}
\times \mathbb{B} \longrightarrow \mathcal{C}$  we
need some combinatorial notions.

\begin{defn}
A digraph $\Gamma$ consists of the following data

\begin{enumerate}
\item A pair of finite sets $(V_{\Gamma},E_{\Gamma})$ called the set of vertices
and edges of $\Gamma$, respectively.
\item A map $(s,t):E_{\Gamma}\longrightarrow V_{\Gamma}\times V_{\Gamma}$.  We call
$s(e)$ and  $t(e)$ the source and target of $e \in V_{\Gamma},$
respectively.
\end{enumerate}

\end{defn}

We use the notations $in(v)=\{e\mid t(e)=v\}$, $i(v)=|in(v)|$,
$out(v)=\{e\mid s(e)=v\}$, and $o(v)=|out(v)|$.  The valence of
$v\in V_{\Gamma}$ is $val(v)=(i(v),o(v))\in
\mathbb{N}^{2}$.  Also we introduce the notation
$V_{\Gamma,in}=\{v\in V_{\Gamma}\mid i(v)=0\}$ and $V_{\Gamma,
out}=\{v\in V_{\Gamma}\mid o(v)=0\}$. Digraphs considered in this
work do not have oriented cycles.  An oriented cycle in $\Gamma$
is a sequence $e_{1},...,e_{n}$ of edges in $\Gamma$ such that
$t(e_{i})=s(e_{i+1})$ for $1\leq i\leq n-1$ and
$t(e_{n})=s(e_{1})$.

\begin{defn}
Let $a$ and $b$ be finite sets. An $(a,b)$-digraph is a triple
$(\Gamma,
\alpha,\beta)$  such that
\begin{enumerate}
\item $\Gamma$ is a digraph.
\item $\alpha: a\longrightarrow V_{\Gamma,in}$ is an injective map.
\item $\beta: b\longrightarrow V_{\Gamma,out}$ is an injective map.
\end{enumerate}
\end{defn}

Let $Digraph(a,b)$ be the groupoid of $(a,b)$-digraphs. A functor
$G:\mathbb{B}^{op}
\times \mathbb{B} \longrightarrow \mathcal{C}$ induces a functor
$G: Digraph(a,b) \longrightarrow  \mathcal{C}$ given by
$$G(\Gamma)=\bigotimes_{v\in V_{\Gamma,
int}}G(in(v),out(v)),$$ $t$ an object of $Digraph(a,b)$ and
$V_{\Gamma, int}=V_{\Gamma}-(\alpha(a)
\sqcup \beta(b))$.

\begin{defn}
The prop $P_{G}$ freely generated by $G:\mathbb{B}^{op}
\times \mathbb{B} \longrightarrow \mathcal{C}$ is given for $n,m \in \mathbb{N}$ by
$$P_{G}(n,m):=
\begin{array}{c}
\underrightarrow{\mathrm{colim}}\,\,G(\Gamma),\\
\end{array}$$

the colimit is taken over the groupoid $Digraph([n],[m]).$
Compositions in $P_G$ are given by gluing digraphs.
\end{defn}
To define props via generators and relations we need to know what
the analogue of an ideal in the prop context is.

\begin{defn}
A subcategory $I$ of $P$ is a prop ideal if
\begin{enumerate}
\item $Ob(I) = Ob(P).$
\item $I(n,m)\otimes P(m,k)\longrightarrow I(n,k)$,
$P(n,m)\otimes I(m,k)\longrightarrow I(n,k)$.
\item $I(n,m)\otimes P(k,l)\longrightarrow I(n\sqcup k,m\sqcup l)$,
$P(n,m)\otimes I(k,l)\longrightarrow I(n\sqcup k,m\sqcup l)$.
\end{enumerate}
for  $n,m,k,l \in \mathbb{N}$.
\end{defn}

We are ready to define Boole as an object in $Prop_{Set}$ . The
prop Boole is a quotient by a prop ideal $I_B$ of the prop freely generated by vertices\\

\begin{figure}[h!]
\begin{center}
\includegraphics[height=1.2cm]{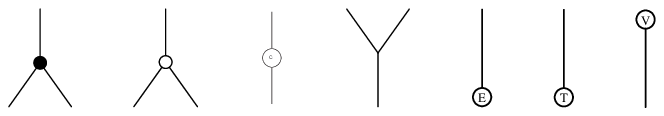}
\end{center}
\end{figure}


representing union, intersection, complement, coproduct, the empty
element, the total element and the valuation, respectively. The
prop ideal $I_{B}$ is generated by the relations given below, each
corresponding with an axiom in the definition of Boolean algebras.

\begin{enumerate}
\item Commutativity for union and intersection

\begin{figure}[h!]
\begin{center}
\includegraphics[height=1.2cm]{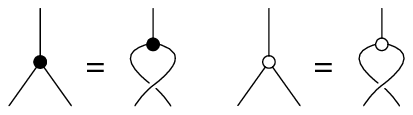}
\end{center}
\end{figure}

\newpage
\item Associativity for union and intersection

\begin{figure}[h!]
\begin{center}
\includegraphics[height=1.2cm]{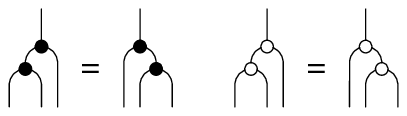}
\end{center}
\end{figure}


\item Distributivity laws

\begin{figure}[h!]
\begin{center}
\includegraphics[height=1.2cm]{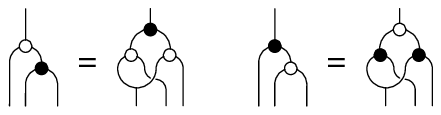}
\end{center}
\end{figure}


\item Properties of the empty and total elements

\begin{figure}[h!]
\begin{center}
\includegraphics[height=1.2cm]{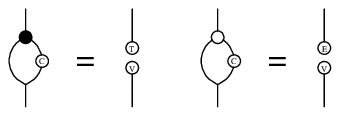}
\end{center}
\end{figure}


\item Absorption Laws

\begin{figure}[h!]
\begin{center}
\includegraphics[height=1.2cm]{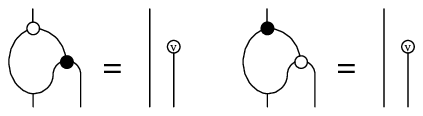}
\end{center}
\end{figure}


\item Coassociativity and cocommutativity

\begin{figure}[h!]
\begin{center}
\includegraphics[height=1.2cm]{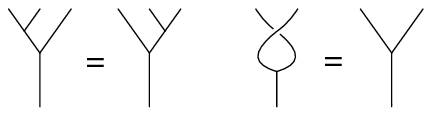}
\end{center}
\end{figure}

\end{enumerate}

\begin{defn}
For $x\in Ob(\mathcal{C})$ we let $End_{x}^{\mathcal{C}}$ be the
prop given by $End_{x}^{\mathcal{C}}(n,m)=\mathcal{C}(x^{\otimes
n},x^{\otimes m}),$ for $n,m \in \mathbb{N}.$
\end{defn}

\begin{defn}
Let $P$ be a $prop$ over $\mathcal{C}$.  A $P$-algebra is a pair
$(x,r)$, where $x$ is an object of $\mathcal{C}$ and
$r:P\longrightarrow End_{x}^{\mathcal{C}}$ is a prop morphism.
\end{defn}

In practice a $P$-algebra  $x$ is given by a family of  maps
$r:P(n,m)\longrightarrow\mathcal{C}(x^{\otimes n},x^{\otimes m})$
satisfying some compatibility conditions.

\begin{thm}
$B$ is a  Boole-algebra  in $Set$ if and only if $x$ is a Boolean
monoid.
\end{thm}
\begin{proof} Assume that $(B,r)$ is a Boole-algebra in $Set$ where
$r: Boole\longrightarrow End_{B}^{Set}$ is a prop morphism. The
images under $r$ of the generators of Boole give operations $\cup,
\cap, (\
\ )^c,t,e,\bigtriangleup, ev$, respectively. For example
$t:\{1\} \longrightarrow B$ and $e:\{1\} \longrightarrow B$ are
identified with elements of $B$. $ev:B \longrightarrow \{1\}$ is
the constant map and plays no essential part in this story.\\

We also get a map $\bigtriangleup: B \longrightarrow B
\times B$ which does seem to fit into the definition of Boolean
monoids. Assume that $\bigtriangleup$ is given by
$\bigtriangleup(a)=(f(a),g(a))$ for $a \in B$.  We use the
relations in $Boole$. The cocommutativity graph implies that
$f=g$. The coassociativity graph implies that $f^{2}=f.$ One of
the absorption graphs implies the identity $f(a) \cup (f(a) \cap
b) = a$ for $a,b \in B$. Thus we obtain

$$f(a)=f^{2}(a) \cup (f^{2}(a) \cap b)=f(a) \cup (f(a) \cap b) = a.$$

Thus $\bigtriangleup(a)=(a,a)$ and it is a simple check that all
other relations in $Boole$ turn $B$ into a Boolen monoid.\\

Assume that $B$ is a Boolean monoid  with operations $\cup,
\cap, (\ \ )^c,$ and distinguished elements $t$ and $e$
that may be thought as maps from $\{1 \}$ to $B$. Take $ev$ to be
the constant map from $B$ to $\{1 \}$, and let $\bigtriangleup$ be
given by $\bigtriangleup(a)=(a,a).$ Let $r$ be the map assigning
to each generator of the Boole prop the corresponding map from the
list above. The fact that $B$ is a Boolean  monoid guarantees that
all the relations defining Boole are satisfied and $r$ extends to
a prop
morphism $r:Boole \longrightarrow End_{B}^{Set}$.\\
\end{proof}

Notice that the functor $<\space{}>: Set \to vect$ induces a
functor $<\space{}>:PROP_{Set}\longrightarrow PROP_{vect}$ given
by $<P>(n,m)=<P(n,m)>$ for  $n,m \in \mathbb{N}$.  The following
result follows from the fact that each generator of the Boole prop
correspond with an operation on Boolean algebras and each relation
in the prop ideal $I_B$ corresponds with an axiom in the
definition of Boolean algebras.

\begin{thm}
$V$ is a $< \mbox{Boole}>$-algebra  in $vect$ if and only if $V$
is a Boolean algebra.
\end{thm}

\section{Symmetric powers of Boolean algebras}\label{spoba}

 The following ideas  introduced in \cite{DiPa} are useful for
studying the symmetric powers of Boolean algebras. Suppose that a
group $G$ acts by automorphisms on the $k$-algebra $A$. The space
of co-invariants $A/G= A/<ga - a \ \ | \ \ g \in G
\mbox { and } a \in A>$ is a $k$-algebra with product given by
$$\overline{a} \overline{b} = \frac{1}{|G|}\sum_{g \in G}\overline{a(gb)}.$$

For each subgroup $K\subset S_{n}$ the Polya functor
$P_{K}:k$-alg$\to k$-alg$ $ from the category of associative
$k$-algebras into itself is defined by: if $A$ is a $k$-algebra
then $P_{K}(A)$ denotes the $k$-algebra whose underlying vector
space is
$$P_{K}(A)=(A^{\otimes n})/\langle a_{1}\otimes\cdot\cdot\cdot\otimes
a_{n}-a_{\sigma^{-1}(1)}\otimes\cdot\cdot\cdot\otimes
a_{\sigma^{-1}(n)}:a_{i}\in A,\sigma\in K\rangle.$$

The rule for the product of $m$ elements in $P_{K}(A)$ is provided
by our next result.

\begin{thm}\label{eddy}
For any $\{a_{ij}\}_{i=1,j=1}^{m,n} \subseteq A$ the following
identity holds in $P_{K}(A)$
$$|K^{m-1}|\prod_{{i=1}}^{m}
\left(\overline{\bigotimes_{{j=1}}^{n}a_{ij}}\right)=\sum_{\sigma\in\{id\}\times
K^{m-1}}\overline{\bigotimes_{{j=1}}^{n}
\left(\prod_{{i=1}}^{m}a_{i\sigma_{{i}}^{-1}(j)}\right)}$$
\end{thm}

In particular for each algebra $A$ and each positive integer $m$
the Polya functor $P_{S_n}$ yields an algebra $P_{S_n}(A)$ which
we denote by $Sym^{n}(A).$ Recall that  $<P[k]>$ denotes the
$k$-vector space generated by the subsets of $[k]$.
 The structural maps $\cup, \cap,$ and $( )^c$ for $<P[k]>$
 are the linear extension of the union, intersection, and complement on $P[k]$.

\begin{defn}
${Sym}^{m}<P[k]>$ is called the symmetric Boolean algebra of type
$(m,k)$. It has operation of union, intersection, and complement
induced by the corresponding operators in $<P[k]>$.
\end{defn}

$S_{x}$ acts by automorphisms on $<P(x)>$ for any finite set $x$.
The next result gives a characterization of the algebra of
co-invariants $<P(x)>/S_{x}.$

\begin{prop}
\begin{enumerate}
\item{$<P(x)>/S_{x}\cong <P[1]>^{\otimes
|x|})/S_{|x|}=Sym^{|x|}<P[1]>$.}
\item{dim$(<P(x)>/S_{x})=|x|+1$.}
\end{enumerate}

\end{prop}

A basis for  $<P[k]>/{S_{k}}$ is given by $\hat{0},...,\hat{k}$
where $\hat{i}$ denotes the equivalence class of $[i] \subseteq
[k]$. Now we study in details the operation of union,
intersection, and complements on the space  $<P[k]>/{S_{k}}$.
Below we use the notation $P(x,k):=\{c\in P(x):|c|=k\}$ for any
set $x$.

\begin{thm} For  $\hat{a},\hat{b}$ in the basis of $<P[k]>/{S_{k}}$
we have
\begin{enumerate}

\item Let $m=min(k-a,b),$
$$\hat{a} \cup \hat{b} =\frac{1}{\left(%
\begin{array}{c}k\\b\end{array}%
\right)}\sum_{l=0}^{m}\left(%
\begin{array}{c}a\\b-l\end{array}%
\right)\left(%
\begin{array}{c} k-a \\ l\end{array}%
\right)\widehat{a+l},$$
\item Let $m=min(a,b),$
$$\hat{a} \cup \hat{b} =\frac{1}{\left(%
\begin{array}{c}k\\b\end{array}%
\right)}\sum_{l=0}^{m}\left(%
\begin{array}{c}a\\l\end{array}%
\right)\widehat{l}.$$
\item $(\hat{a})^{c}=\widehat{k-a}.$
\end{enumerate}
\end{thm}

\begin{proof}

\begin{enumerate}

\item \begin{eqnarray*}
\hat{a} \cup \hat{b}&=&\frac{1}{k!}\sum_{\sigma\in S_{k}}\overline{[a]\cup\sigma[b]}=\frac{1}{\left(%
\begin{array}{c}k \\ b\end{array}%
 \right)}\sum_{c\in P([k],b)}\overline{[a]\cup c} =\\
 &=&\frac{1}{\left(%
 \begin{array}{c}k \\ b\end{array}%
 \right)}\sum_{\begin{subarray}{1}c_{0}\subset P([k]-[a],l)\\c_{1}\subset P([a],b-l)\end{subarray}}\overline{[a]\cup c_{0}}=\frac{1}{\left(%
 \begin{array}{c}k \\ b \end{array}%
 \right)}\sum_{{l=0}}^{m}{\left(%
 \begin{array}{c}a \\ b-l \end{array}%
 \right)}{\left(%
 \begin{array}{c}k-a \\ l \end{array}%
 \right)}\widehat{a+l}.\\
\end{eqnarray*}

\item
Follows from the fact that the number of permutations $ \sigma \in
S_k$ such $|[a] \cap \sigma([b])|=l$ is given by

$$ \left( \begin{array}{c} a \\ l\end{array}\right)
 \left( \begin{array}{c} b \\ l\end{array}\right) l!
\left( \begin{array}{c} k-a \\ b-l\end{array}\right)(b-l)!(k-b)!$$
\item Obvious.

\end{enumerate}
\end{proof}

Let $\pi= \{b_1,...,b_k \}$ be a partition of $x$ and $S_{\pi}
\subseteq S_x$ the Young
subgroup consisting of block preserving permutations of $x$. Our
next result characterizes algebras of the form $<P(x)>/S_{\pi}$.

\begin{prop}
\begin{enumerate}
\item{$<P(x)>/S_{\pi}\cong \bigotimes_{i=1}^{k}<P[1]>^{\otimes
|b_i|}/S_{|b_i|}=\bigotimes_{i=1}^{k}Sym^{|b_i|}<P[1]>$.}
\item{dim$(<P(x)>/S_{\pi})=\prod_{i=1}^{k} (|b_i|+1)$.}
\end{enumerate}
\end{prop}

We close this section by taking a closer look at the symmetric
Boolean algebra $Sym^{2}<P[1]>$ and the cyclic Boolean algebra
$<P[1]>^{\otimes 3}/\mathbb{Z}_3$. \\

The space $Sym^{2}<P[1]>$ has basis
$\widehat{0}=\overline{(\emptyset,
\emptyset)},\widehat{1}=\overline{([1], \emptyset)}$ and $\widehat{2}=\overline{([1], [1])}$.
The union $\cup: Sym^{2}<P[1]> \otimes Sym^{2}<P[1]>$
$\longrightarrow Sym^{2}<P[1]>$ is given for $i=0,1,2$ by
\begin{itemize}
\item{$\widehat{0} \cup \widehat{i} = \widehat{i}$.}
\item{$\widehat{1} \cup \widehat{1} = \frac{1}{2}\widehat{1} + \frac{1}{2}\widehat{2}$.}
\item{$\widehat{2} \cup \widehat{i} = \widehat{2}$.}
\end{itemize}
The intersection  $\cap: Sym^{2}<P[1]> \otimes Sym^{2}<P[1]>$
$\longrightarrow Sym^{2}<P[1]>$ is given for $i=0,1,2$ by
\begin{itemize}
\item{$\widehat{0} \cap \widehat{i} = \widehat{0}$.}
\item{$\widehat{1} \cap \widehat{1} = \frac{1}{2}\widehat{0} + \frac{1}{2}\widehat{1}$.}
\item{$\widehat{2} \cap \widehat{i} = \widehat{i}$.}
\end{itemize}

The complement $(\ \ )^{c}: Sym^{2}<P[1]> \longrightarrow
Sym^{2}<P[1]>$ is given by
\begin{itemize}
\item{$\widehat{0}^c= \widehat{2}$, $\widehat{1}^{c} = \widehat{1}$, and $\widehat{2}^c= \widehat{0}$.}
\end{itemize}

Although the algebra $Sym^{2}([1])$ does not satisfy all  axioms
required to make it into a Boolean algebra (absorption fails) it
does share many of the properties of Boolean algebras,
and in any event it is a mathematical object of great interest.\\

Let us consider in details the  third cyclic power of the Boolean
algebra $<P[1]>$, namely $<P[1]>^{\otimes 3}/\mathbb{Z}_3$. It has
basis $\widehat{0}=\overline{(\emptyset,\emptyset,\emptyset)}$,
$\widehat{1}=\overline{([1],\emptyset,\emptyset)}$,
$\widehat{2}=\overline{([1],[1],\emptyset)}$ and
$\widehat{3}=\overline{([1],[1],[1])}.$ The union $\cup:
<P[1]>^{\otimes 3}/\mathbb{Z}_3 \otimes <P[1]>^{\otimes
3}/\mathbb{Z}_3
\longrightarrow P[1]^{\otimes 3}/\mathbb{Z}_3$ is given for
$i=0,1,2,3$ by
\begin{itemize}
\item{$\widehat{0} \cup \widehat{i} = \widehat{i}$.}
\item{$\widehat{1} \cup \widehat{1} = \frac{1}{3}\widehat{1} + \frac{2}{3}\widehat{2}.$}
\item{$\widehat{1} \cup \widehat{2} = \frac{2}{3}\widehat{2} + \frac{1}{3}\widehat{3}.$}
\item{$\widehat{2} \cup \widehat{2} = \frac{1}{3}\widehat{2} + \frac{2}{3}\widehat{3}.$}
\item{$\widehat{3} \cup \widehat{i} = \widehat{3}$.}
\end{itemize}

The intersection $\cap: <P[1]>^{\otimes 3}/\mathbb{Z}_3 \otimes
<P[1]>^{\otimes 3}/\mathbb{Z}_3
\longrightarrow <P[1]>^{\otimes 3}/\mathbb{Z}_3$ is given for
$i=0,1,2,3$ by
\begin{itemize}
\item{$\widehat{0} \cap \widehat{i} = \widehat{0}$.}
\item{$\widehat{1} \cap \widehat{1} = \frac{2}{3}\widehat{0} + \frac{1}{3}\widehat{1}.$}
\item{$\widehat{1} \cap \widehat{2} = \frac{1}{3}\widehat{0} + \frac{2}{3}\widehat{1}.$}
\item{$\widehat{1} \cap \widehat{2} = \frac{2}{3}\widehat{1} + \frac{1}{3}\widehat{2}.$}
\item{$\widehat{3} \cap \widehat{i} = \widehat{i}$.}
\end{itemize}

The complement map $(\ \ )^{c}: P([1])^{\otimes 3}/\mathbb{Z}_3
\longrightarrow P([1])^{\otimes 3}/\mathbb{Z}_3$ is given by
\begin{itemize}
\item{$[0]^c= [3]$, $[1]^c= [2]$, $[2]^c= [1]$, and $[3]^c= [0]$ .}
\end{itemize}

\section{Symmetric inclusion-exclusion principles}

In this Section  we take $k= \mathbb{R}$. We write
$\{a_{1},...,a_{m}\}$ for the basis element
$\overline{a_{1}\otimes\cdot\cdot\cdot\otimes a_{m}}\in
<P[k]>^{\otimes m}/S_{m}$. The following result follows from
Theorem \ref{eddy}.\\

\begin{thm}\label{verde}
Let $\{a_{1}^i,...,a_{m}^i\}$ be in the basis of $\in
<P[k]>^{\otimes m}/S_{m}$ for $1\leq i\leq n$. The union in
$<P[k]>^{\otimes m}/S_{m}$ is given by
$$\bigcup_{{i=1}}^n
\{a_{1}^i,...,a_{m}^i\}=\frac{1}{(m!)^{n-1}}\sum_{\begin{subarray}{1}\sigma\in\{1\}\times
S_{m}^{(n-1)}\end{subarray}}\{\bigcup_{{i=1}}^{n}a_{{{\sigma_{i(1)}}}}^{i},...,
\bigcup_{{i=1}}^{n}a_{{{\sigma_{i(m)}}}}^{i}\}$$
\end{thm}

\begin{exmp}

For $m,n=2$ we get
$$\{a_{1}^{1},a_{2}^{1}\}\cup\{a_{1}^{2},a_{2}^{2}\}=\frac{1}{2}\{a_{1}^{1}\cup a_{1}^{2},a_{2}^{1}\cup a_{2}^{2}\}
+\frac{1}{2}\{a_{1}^{1}\cup a_{2}^{2},a_{2}^{1}\cup a_{1}^{2}\}.$$
In a better notation
$$\{a,b\}\cup\{c,d\}=\frac{1}{2}\{a\cup c,b\cup d\}
+\frac{1}{2}\{a\cup d,b\cup c\}.$$
\end{exmp}

A measure on a finite set $x$ is a map {$\mu:P(x)\rightarrow
\mathbb{R}$} such that $\mu(a \cup b)=\mu(a)+\mu(b)$ for $a,b
 \subseteq x$ disjoint. Fix a measure $\mu$ on $[k]$. An element  $\{a_{1},...,a_{m}\}$
 in the basis of $<P[k]>^{\otimes m}/S_{m}$ determines a vector
$(\mu(a_{1}),...,\mu(a_{n}))\in\mathbb{R}^{m}/S_{m}$.  Functions
on $\mathbb{R}^{m}/S_{m}$ are known as symmetric functions. There
are many interesting examples of polynomial symmetric functions
such as the power functions, the elementary symmetric functions,
the homogeneous functions, the Schur functions and so on. For
example the polynomial $x_{{1}}^{l}+\cdot\cdot\cdot+x_{{m}}^{l}$
is $S_{m}$-invariant. Each symmetric function can be used to
obtain a symmetric form of the inclusion-exclusion principle. The
reader will find interesting information on the
inclusion-exclusion principle and its generalizations in several
papers by Rota and his collaborators in \cite{rota}. We use the
inclusion-exclusion principle in the following form.

\begin{prop}\label{bit}
Let $a_{1},...,a_{n}\in P(x)$ then
$|\bigcup_{i=1}^{n}a_{i}|=\sum_{I\subseteq[n]}(-1)^{|I|+1}|\bigcap_{i\in
I}a_{i}|.$
\end{prop}
We consider the symmetric inclusion-exclusion principles
\footnote{ In \cite{Ge} Gessel uses the name symmetric inclusion-exclusion to
refer to a different mathematical gadget.}  derived from the
power, elementary, and homogeneous symmetric functions. Other
symmetric functions can be used as well but we shall not do so
here.  The power function $p_{l}:<P[k]>^{\otimes m}/S_{m}
\rightarrow
\mathbb{R}$ is given on the basis by
$p_{l}(\{a_{1},...,a_{m}\})=\sum_{{i=1}}^{m}\mu(a_{i})^{l}.$ We
use the power functions  $p_l$  to get a symmetric form of the
inclusion-exclusion principle.

\begin{thm}\label{siep1}
Let $\{a_{1}^i,...,a_{m}^i\}$ be in the basis of $<P[k]>^{\otimes
m}/S_{m}$ for $1\leq i\leq n$. Then
$$p_{l}(\bigcup_{{i=1}}^n \{a_{1}^i,...,a_{m}^i\})=\frac{1}{(m!)^{n-1}}\
sum_{\begin{subarray}{1}\sigma\in\{1\}\times S_{m}^{(n-1)}\\ j\in \{1,...,m\}\\ \Sigma
c_I=l\\ \end{subarray}}\left(%
\begin{array}{c}l \\ \{c_I\}\end{array}%
\right)\prod_{I\subseteq [n]}(-1)^{(\mid I\mid +1)c_I}\mu(\bigcap_{i\in I}
a_{\sigma_{i}(j)}^{i})^{c_I}.$$
\end{thm}

\begin{proof}
\begin{eqnarray*}
p_{l}(\bigcup_{{i=1}}^n\{a_{1}^i,...,a_{m}^i\})&=&\frac{1}{(m!)^{n-1}}\sum_{\begin{subarray}{1}\sigma\in\{1\}\times
S_{m}^{(n-1)}\\ j\in \{1,...,m\}
\end{subarray}}\mu(\bigcup_{{i=1}}^{n}a_{{{\sigma_{i(j)}}}}^{i})^{l}\\
&=&\frac{1}{(m!)^{n-1}}\sum_{\begin{subarray}{1}\sigma\in\{1\}\times
S_{m}^{(n-1)}\\ j\in \{1,...,m\}
\end{subarray}}\left(\sum_{\begin{subarray}{1}I\subseteq
[n]\end{subarray}}(-1)^{|I|+1}\mu(\bigcap_{i\in I}a_{{\sigma_{i(j)}}}^{i})\right)^{l}\\
&=&\frac{1}{(m!)^{n-1}}\sum_{\begin{subarray}{1}\sigma\in\{1\}\times
S_{m}^{(n-1)}\\ j\in \{1,...,m\}
\end{subarray}}\left(\sum_{\begin{subarray}{1}\Sigma c_{I}=l\\
\end{subarray}}\left(\begin{array}{c}l \\\{c_{I}\}\end{array}\right)\prod_{I\subseteq [n]}[(-1)^{|I|+1}\mu(\bigcap_{i\in
I}a_{{\sigma_{i(j)}}}^{i})]^{c_{I}}\right)\\
&=&\frac{1}{(m!)^{n-1}}\sum_{\begin{subarray}{1}\sigma\in\{1\}\times
S_{m}^{(n-1)}\\ j\in \{1,...,m\}\\ \Sigma
c_I=l\\ \end{subarray}}\left(%
\begin{array}{c}l \\ \{c_I\}\end{array}%
\right)\prod_{I\subseteq [n]}(-1)^{(\mid I\mid +1)c_I}\mu(\bigcap_{i\in I}
a_{\sigma_{i}(j)}^{i})^{c_I}.\\
\end{eqnarray*}
\end{proof}

\begin{cor}
For $l=1$ we have
$$p_1(\bigcup_{{i=1}}^n \{a_{1}^i,...,a_{m}^i\})=\frac{1}{(m!)^{n-1}}\sum_{\begin{subarray}{1}\sigma\in\{1\}\times
S_{m}^{(n-1)}\\ j\in \{1,...,m\}\\I\subseteq [n] \end{subarray}}
(-1)^{\mid I\mid +1}\mu(\bigcap_{i\in I}a_{\sigma_{i}(j)}^{i}).$$
\end{cor}

\begin{cor}
For $l=1$, $n=2$ we have
$$p_1(\{a_{1}^1,...,a_{m}^1\} \cup \{a_{1}^1,...,a_{m}^2\})=\frac{1}{m!}\sum_{\begin{subarray}{1}\sigma\in S_{m}\\ j\in [m]\\ \end{subarray}}
\{\mu(a_{j}^{1}) +  \mu(a_{\sigma(j)}^{2}) - \mu(a_{j}^{1}\bigcap a_{\sigma_(j)}^{2})\}.$$
\end{cor}

\smallskip

A generalized inclusion-exclusion principle using the elementary
symmetric functions
\begin{eqnarray*}
e_{l}(x_1,...,x_m)&=&\sum_{1\leq t_{1} < t_{2} <
\cdot\cdot\cdot <
            t_{l}\leq m}\prod_{j=1}^{l}x_{t_{j}}.\\
\end{eqnarray*}

is given by

\begin{thm}\label{siep}
Let $\{a_{1}^i,...,a_{m}^i\}$ be in the basis of $<P[k]>^{\otimes
m}/S_{m}$ for $1\leq i\leq n$. Then
$$e_{l}(\bigcup_{{i=1}}^n \{a_{1}^i,...,a_{m}^i\})=\frac{1}{(m!)^{n-1}}\sum_{\begin{subarray}{c}
\sigma\in \{1\}\times S_{m}^{n-1} \\1\leq t_{1}<t_{2}<\cdots<t_{l}\leq m\\ f:[l]\to P([n])
\end{subarray}}\prod_{j=1}^l(-1)^{\mid f(j)\mid
+1}\mu(\bigcap_{i\in f(j)}a_{\sigma_{i}(t_{j})}^i).$$
\end{thm}

\begin{proof} According to Theorem \ref{verde} we get

\begin{eqnarray*}
e_{k}(\bigcup_{i=1}^n\{a_{1}^i,...,a_{m}^i\})&=&\frac{1}{(m!)^{n-1}}\sum_{\begin{subarray}{c}\sigma\in
\{1\}\times S_{m}^{n-1} \\1\leq p_{1}<p_{2}<\cdots<p_{l}\leq
m\end{subarray}}\prod_{j=1}^l\mu(\bigcup_{i\in I}a_{\sigma_{i}(p_j)}^i)\\
&=&\frac{1}{(m!)^{n-1}}\sum_{\begin{subarray}{c}\sigma\in
\{1\}\times S_{m}^{n-1} \\1\leq p_{1}<p_{2}<\cdots<p_{l}\leq
m\end{subarray}}\prod_{j=1}^l\sum_{I\subseteq [n]}(-1)^{\mid I\mid
+1}\mu(\bigcap_{i\in I}a_{\sigma_{i}(p_j)}^i)\\
&=&\frac{1}{(m!)^{n-1}}\sum_{\begin{subarray}{c}\sigma\in
\{1\}\times S_{m}^{n-1}\\1\leq p_{1}<p_{2}<\cdots<p_{l}\leq m\\
f:[l]\to P([n])
\end{subarray}}\prod_{j=1}^l(-1)^{\mid f(j)\mid
+1}\mu(\bigcap_{i\in f(j)}a_{\sigma_{i}(p_{j})}^i).\\
\end{eqnarray*}
\end{proof}
\smallskip

\begin{exmp}
Let  $n=2$, $m=2$ and $l=2$. The map $e_2:<P[k]>^{\otimes 2}/S_{2}
\longrightarrow \mathbb{R}$ is given by
$e_2(\{a,b\})=\mu(a)\mu(b)$ for $a,b \in P[k]$.  Theorem
\ref{siep} implies that

\begin{eqnarray*}
2e_{2}(\{a,b\}\cup \{c,d\})&=&2\mu(a)\mu(b)+2\mu(c)\mu(d)+\mu(a)\mu(d)+\mu(c)\mu(b)\\
& &+\ \ \mu(a)\mu(c)+\mu(d)\mu(b)-\mu(a)\mu(b\cap d)+\mu(c)\mu(b\cap d)\\
& &+\ \ \mu(b)\mu(a\cap c)+\mu(d)\mu(a\cap c)+\mu(a)\mu(b\cap c)\\
& &+\ \ \mu(d)\mu(b\cap c)+\mu(b)\mu(a\cap d)+\mu(c)\mu(a\cap d).
\end{eqnarray*}

\end{exmp}

The generalization of the inclusion-exclusion principle using the
homogenous symmetric functions
\begin{eqnarray*}
h_{l}(x_1,...,x_m)&=&\sum_{1\leq t_{1} \leq t_{2} \leq
\cdot\cdot\cdot <t_{l}\leq m}\prod_{j=1}^{l}x_{t_{j}}.\\
\end{eqnarray*}
is given by

\begin{thm}\label{siep}
Let $\{a_{1}^i,...,a_{m}^i\}$ be in the basis $<P[k]>^{\otimes
m}/S_{m}$ for $1\leq i\leq n$. Then
$$h_{l}(\bigcup_{{i=1}}^n \{a_{1}^i,...,a_{m}^i\})=\frac{1}{(m!)^{n-1}}\sum_{\begin{subarray}{c}
\sigma\in \{1\}\times S_{m}^{n-1} \\1 \leq t_{1} \leq t_{2} \leq \cdots \leq t_{l}\leq m\\ f:[l]\to P([n])
\end{subarray}}\prod_{j=1}^l(-1)^{\mid f(j)\mid
+1}\mu(\bigcap_{i\in f(j)}a_{\sigma_{i}(t_{j})}^i).$$
\end{thm}

\begin{exmp}
Let  $n=2$, $m=2$ and $l=2$. The map $h_2:<P([k])>^{2}/S_{2}
\longrightarrow \mathbb{R}$ is given by $h_2(\{a,b\})=
\mu(a)^{2} + \mu(a)\mu(b) + \mu(b)^{2}$ for $a,b
\in P[k]$. Theorem
\ref{siep} implies that

\begin{eqnarray*}
2h_{2}(\{a,b\}\cup \{c,d\})&=&[\mu(a)+\mu(c)-\mu(a\cap c)]^{2}+[\mu(b)+\mu(d)-\mu(b\cap d)]^{2}\\
& &+\ \ [\mu(a)+\mu(d)-\mu(a\cap d)]^{2}+[\mu(b)+\mu(c)-\mu(b\cap c)]^{2}\\
& &+\ \ 2\mu(a)\mu(b)+2\mu(c)\mu(a)+\mu(a)\mu(d)+\mu(c)\mu(b)+\mu(a)\mu(c)\\
& &+\ \ \mu(d)\mu(a)-\mu(a)\mu(b\cap d)+\mu(c)\mu(b\cap d)+\mu(b)\mu(a\cap c)\\
& &+\ \ \mu(d)\mu(a\cap c)+\mu(a)\mu(b\cap c)+\mu(d)\mu(b\cap c)+\mu(b)\mu(a\cap d)\\
& &+\ \ \mu(c)\mu(a\cap d).
\end{eqnarray*}

\end{exmp}

\section{Propositional logic and Boolean algebras}

It is hard to do any work on Boolean monoids and not to mention at
all its relation with propositional logic. Indeed the motivation
of Boole himself to introduce Boolean monoids was to describe the
mathematical structures that control the laws of though.
Propositional logic deals with the relation of deduction among
sequences of sets of sentences constructed from a given finite set
of propositions connected by a fixed set of connecting symbols.
Let us denote the set of given propositions $C$ and the set of
sentences by $S$. There are many ways to describe a system of
propositional logic but in any of them one can imagine that there
exists a sort of logical agent capable of performing the following
tasks

\begin{itemize}
\item{Recognize when a grammatical construction is  an element of $S$.
The agent is able to translate into sentences in $S$ expressions
of the form $s\vee t$, $s\wedge t$, and $-s$ for sentences $s$ and
$t$ in $S$.}

\item{Decide wether or not a sequence of sets of sentences $c_1,...,c_n$
is a deduction. }

\item{Assign a truth value to each sentence in $S$ when provided with a assignment
of truth values for propositions in $P$, i.e., an element of $\{0,1\}^C$.}\\

A sentence $s$ is said to imply a sentence $t$ if there exists a
deduction $c_1,...,c_n$ such that $c_1 = \{ s\}$ and $c_n = \{
t\}$. The logical operator is said to be sound and complete if the
following property holds
\item{Sentence $s \in S$ implies sentence $t \in S$ if for any assignment of truth values to propositions in
 $C$  the truth value of $t$ is $1$ if the truth value of $s$ is $1$. It is no hard to show the existence of sound a complete logical
agents \cite{smu}.}
\end{itemize}
Boolean monoids appear within the context of propositional logic
as follows. Say that sentences $s$ and $t$ in $S$ are equivalent
if $s$ implies $t$ and $t$ implies $s$. Let $B(S)$ be the quotient
of $S$ by this equivalence relation. $B(S)$ comes equipped with a
natural structure of Boolean monoid with operations defined by
$[s]
\cup[t]=[s \vee t]$, $[s] \cup[t]=[s \wedge t]$, and $[s]^c =
[-s]$, for $[s]$ and $[t]$ in $B(S)$ . The total element is $[s
\vee -s]$ and the empty element is $[s \wedge -s].$  The Boolean
monoid $B(S)$ is isomorphic to the Boolean monoid $P(\{0,1\}^C)$
via the map
$$m:B(S)\rightarrow P(\{0,1\}^C)$$ sending each sentence $[s] \in S$ into the set of its models
$$m([s])=\{v \in \{0,1\}^C \ \ | \ \ \mbox{ the truth value of $s$ according to $v$ is $1$
}\}.$$

Summarizing sentences in $S$ describe subsets of $\{0,1\}^C$ and
two sentences describe the same set if and only if they are
equivalent. The power of the logical description of $P(\{0,1\}^C)$
lies in the possibility of describing the same set in a variety of
ways. For example the logical agent may be told that a subset of
$\{0,1\}^C$ is described by a sentence $s$, that another subset of
$\{0,1\}^C$ is described by a sentence $t$, and be asked to
provide a sentence which describes the union of those sets.  It
will readily
answer that $s \vee t$ is the sought after sentence.\\

It is natural to wonder if any logical meaning can be ascribed to
the Boolean algebra $<B(S)>$. Although preliminary we venture an
answer: assume the logical agent is told that a sentence $s_i$
describes an unknown subset of $\{0,1\}^C$ with probability $p_i$
for $1 \leq i \leq n$, and that a sentence $t_j$ describes another
unknown subset of $\{0,1\}^C$ with probability $q_j$ for
 $1 \leq j \leq m.$ If asked to find a sentence that describes the
 union of those subsets the logical agent will answer: the
 sentence $s_i \vee t_j$  describes the union of the unknown sets
 with probability $p_{i}q_{j}.$ This is the only consistent answer
 with the product rules on $<B(S)>$ which is given by
 $$(\sum_{i=1}^{n}p_{i}[s_i])\cup (\sum_{j=1}^{m}q_{j}[t_j])
 = \sum_{i=1,j=1}^{n,m}p_{i}q_{j}[s_i \vee t_j].$$

This probabilistic interpretation  applies as well to the Boolean
algebra $<P(x)>$. Let $v$ and $w$ be a couple of vectors in
$<P(x)>$ given by $v=
\sum_{a \subseteq x}v_{a}a \mbox{ and } w= \sum_{b
\subseteq x}v_{b}b.$ Assume that the coefficients of $v$ and
$w$, respectively, are positive and add to one. This allow us to
think that $v_a$ represents the probability that the unknown
subset $v$ of $x$ be equal to $a$. Similarly $w_b$ represents the
probability that $w$ be equal to $b$. Under this conditions we
have that

\begin{itemize}
\item{The probability that $v
\cup w$ be equal to $c$ is given by $(v \cup w)_c = \sum_{a \cup b
=c}v_{a}w_{b}$.}

\item{
The probability that $v
\cap w$ be equal to $c$ is given by $(v \cap w)_c = \sum_{a \cap b
=c}v_{a}w_{b}$.}

\item{The probability that $v^{c}$ be equal to $a$ is $v_{a^c}$.}
\end{itemize}

Finally we invite the reader to take another look at the
structural coefficients of the algebras $Sym^{2}<P[1]>$ and
$<P[1]>^{\otimes 3}/\mathbb{Z}_3$ given in Section \ref{spoba} and
check that they are indeed consistent with the probabilistic
interpretation just outlined.

\subsection*{Acknowledgment} Thanks to Mauricio Angel, H\'ector Blandin,
Edmundo Castillo and Eddy Pariguan. This paper is dedicated to the
memory of Professor Ramon Castillo Ariza.

\smallskip

\noindent Rafael D\'\i az\\
ragadiaz@gmail.com\\

\noindent Mariolys Rivas\\
mariolysrivas07@gmail.com \\
Universidad Central de Venezuela\\

\end{document}